\documentclass[12pt]{article}
\usepackage{latexsym}
\usepackage{amssymb}
\usepackage{amsthm}
\usepackage{amsmath}
\usepackage{graphicx}
\usepackage{comment}
\def\calp{{\cal P}}
\def\cald{{\cal D}}
\newcommand{\lang}{\langle}
\newcommand{\rang}{\rangle}
\newcommand{\Q}{\mathbb{Q}}

\newtheorem{thm}{Theorem}[section]
\newtheorem{theorem}{Theorem}[section]
\newtheorem{lemma}[thm]{Lemma}

\newtheorem{conjecture}[thm]{Conjecture}
\newtheorem{corollary}[thm]{Corollary}

\begin{document}

\renewcommand{\theequation}{\arabic{section}.\arabic{equation}}
\thispagestyle{empty}

\vskip 20pt
\begin{center}
{\bf MINIMAL BAR TABLEAUX}
\vskip 15pt
{\bf Peter Clifford\footnote{
Homepage:  {\tt http://cnri.dit.ie/}}
}\\
{\it CNRI,}
{\it Dublin Institute of Technology,
Ireland}\\
{\tt peter@cnri.dit.ie}
\vskip 10pt
\end{center}

\begin{abstract} 
Motivated by Stanley's results in \cite{St02}, we generalize the rank of a partition $\lambda$ to the rank of a
shifted partition $S(\lambda)$.  We show that the number of bars required in
a minimal bar tableau of $S(\lambda)$ is max$(o, e + ( \ell(\lambda)
\ \mathrm{mod} \   2))$,
where $o$ and $e$ are the number of odd and even rows of $\lambda$.  As a   
consequence we show that the irreducible projective characters of
$S_n$ vanish on certain conjugacy classes.  Another corollary
is a lower
bound on the degree of the terms in the expansion of Schur's $Q_{\lambda}$
symmetric functions in terms of the power sum symmetric functions.
\end{abstract}

\section{Introduction}

Let $\lambda = (\lambda_1, \lambda_2, \ldots)$ be a partition of the
integer $n$, i.e., $ \lambda_1 \ge \lambda_2 \ge \cdots \ge 0$ and
$\sum \lambda_i =n$.   
The
\emph{length} $\ell(\lambda)$ of a partition $\lambda$ is the number of 
nonzero parts of $\lambda$.  
The (Durfee or Frobenius) \emph{rank} of
$\lambda$, denoted rank($\lambda$), is the length of the main diagonal
of the diagram of $\lambda$, or equivalently, the largest integer $i$
for which $\lambda_i \ge i$.  The rank of $\lambda$ is the least
integer  $r$ such that
$\lambda$ is a disjoint union of $r$ border strips (also called 
ribbons or rim hooks).  
Denote by $\chi ^ \lambda (\pi)$ the character of the irreducible
representation of $S_n$ indexed by $\lambda$ evaluated at a
permutation of cycle type $\pi$.
An easy application of the Murnaghan-Nakayama rule shows that if $\ell(\pi)<
\mbox{rank}(\lambda)$, we must have $\chi ^ \lambda (\pi) = 0$.  
As a corollary, the expansion of the Schur function in terms of the
power sum symmetric functions
$$s_{\lambda} = \sum_{\pi} \chi^{\lambda} (\pi)
\frac{p_{\pi}}{z_{\pi}}$$ 
contains only terms 
$\chi^{\lambda} (\pi)\frac{
p_{\pi}}{z_{\pi}}$ such that $\ell(\pi) \ge $ rank $(\lambda)$.

We generalize the notion of rank to shifted diagrams
$S(\lambda)$ of a partition with distinct parts by counting the minimal number
of bars in a bar tableau.
Using Morris' projective analogue of the Murnaghan-Nakayama rule,  we
show that the
irreducible projective characters of $S_n$ vanish on conjugacy classes
indexed by partitions with few parts.
This enables us  to give a lower bound on the length of the $\mu$
which appear in the expansion of the Schur
Q-functions in terms of the $p_\mu$.


\section{Bar Tableaux}

\begin{figure}[h]
\centering
\includegraphics[width=2.0in]{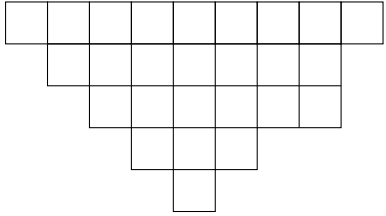}
\caption{The shifted diagram of the shape 97631}
\label{fig:shiftedshape}
\end{figure}

Let $\cald (n)$ be the set of all partitions of $n$ into distinct parts.
The \emph{shifted diagram}, $S(\lambda)$, of shape $\lambda$ is obtained by
forming $l$ rows of nodes with $\lambda_i$ nodes in the $i$th row
such that, for all $i>1$, the first node in row $i$ is placed
underneath the second node in row $(i-1)$.  For instance Figure
\ref{fig:shiftedshape} shows the shifted diagram of the shape 97631.

We follow the treatment of Hoffman and Humphreys \cite{Ho92} to define
bar tableaux.  These occur in the inductive formula for the
projective characters of $S_n$, first proved by Morris
\cite{Mo62}.  Let $r$ be 
an odd positive integer, and let 
$\lambda \in \cald(n)$ 
have length $l$.  Below we define:
\begin{enumerate}
\item a subset, $I_+ \cup I_0 \cup I_- = I(\lambda,r)$, of integers
  between $1$ and $l$; and
\item for each $i \in I(\lambda,r)$, a strict partition $\lambda(i,r)$
  in $ \cald(n-r)$ (despite the notation, $\lambda(i,r)$ is a
  function of $\lambda$, as well as of $(i,r)$). 
\end{enumerate}
Let
$$I_+ = \{i:\lambda_{j+1} < \lambda_i -r < \lambda_j \mbox{ for some }
j \le l, \mbox{ taking } \lambda_{l+1}=0 \}.$$
In other words $I_+$ is the set of all rows of $\lambda$ which we can
remove $r$ squares from and still leave a composition with distinct
parts.  For example, if $r=5$ and $\lambda=97631$, then $I_+ =
\{1,2\}$.  If $i \in I_+$, then $\lambda_i > r$, and we define
$\lambda(i,r)$ to 
be the partition obtained from $\lambda$ by removing $\lambda_i$ and
inserting $\lambda_i -r$ between $\lambda_j$ and
$\lambda_{j+1}$. Continuing our example above, $\lambda(2,5)=96321$.  Let
$$ I_0 = \{i:\lambda_i = r \},$$
which is empty or a singleton.  For $i \in I_0$, remove $\lambda_i$
from $\lambda$ to obain $\lambda(i,r)$. Let
$$ I_- = \{ i: r-\lambda_i = \lambda_j \mbox{ for some } j \mbox{ with
  } i<j  \le l \}.$$
Equivalently $I_-$ is the set of all rows of $\lambda$ for which
there is some shorter row of $\lambda$ such that the total number of
squares in both rows is $r$.  For example, if $r=7$ and
$\lambda=97631$, then  $I_- =
\{3\}$.  If $i \in I_-$, then
$\lambda_i<r$, and $\lambda(i,r)$ is formed by 
removing both $\lambda_i$ and $\lambda_j$ from $\lambda$.

\begin{figure}[h]
\centering
\includegraphics[width=2.0in]{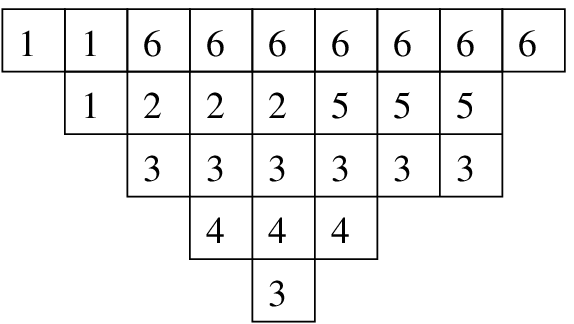}
\caption{A bar tableau of the shape 97631}
\label{fig:bartableau}
\end{figure}

For each $i \in  I(\lambda,r)$ the associated \emph{r-bar} is given as
follows.  If $i$ is in $I_+$ or $I_0$, the $r$-bar consists of the
rightmost $r$ nodes in the $i$th row of $S(\lambda)$.  We say the
$r$-bar is of \emph{type 1} or \emph{type 2} respectively.  For
example, the squares in Figure \ref{fig:bartableau} labelled by 6 are a
$7$-bar of type 1.  The squares labelled by 4 are a $3$-bar of type 2.
If $i$ is in $I_-$,  the $r$-bar
consists of all the nodes in both the $i$th and $j$th rows, a total of
$r$ nodes.  We say the $r$-bar is of \emph{type 3}.  The squares in
Figure \ref{fig:bartableau} labelled by 3 are a $7$-bar of type 3.


Define a \emph{bar tableau} of shape $\lambda$ to be an assignment of
positive integers to the squares of $S(\lambda)$ such that
\begin{enumerate}
\item the set of squares occupied by the biggest integer is an
  $r$-bar $B$, and
\item if we remove the $r$-bar $B$ and reorder the rows, the result is
  a bar  tableau.
\end{enumerate}

Equivalently we can define a \emph{bar tableau} of shape $\lambda$ to
be an assignment  of
positive integers to the squares of $S(\lambda)$ such that
\begin{enumerate}
\item the entries are weakly increasing across rows,
\item each integer $i$ appears an odd number of times,
\item $i$ can appear in at most two rows; if it does, it must
  begin both rows (equivalent to the bar being of type 3),
\item the composition remaining if we remove all squares
  labelled by integers larger than some $i$ has distinct parts.
\end{enumerate}

\begin{figure}[h]
\centering
\includegraphics[width=5.0in]{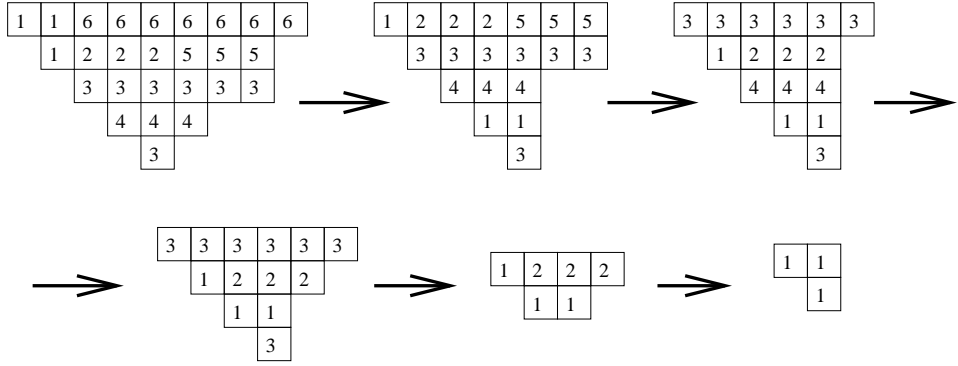}
\caption{Checking legality of the bar tableau of the shape 97631}
\label{fig:tableauchain}
\end{figure}

For example, Figure \ref{fig:tableauchain} shows the chain of
partitions remaining if we remove all squares labelled by integers
larger than some $i$ from the tableau in Figure
\ref{fig:bartableau}.  This demonstrates the legality of that tableau.

\section{Minimal Bar Tableaux}

We introduce an operation on minimal bar
tableaux which preserves the number of bars, and prove some facts
about tableaux resulting from this operation.  A bar tableau of
$\lambda$ is \emph{minimal} if the number of bars is minimized,
i.e. there does not exist a bar tableau with fewer bars.

\begin{lemma}\label{l:noevenbdr}
There exists a minimal bar tableau $T^*$ such that there is no bar
boundary an even number of squares along any row.
\end{lemma}

\begin{figure}[h]
\centering
\includegraphics[width=4.0in]{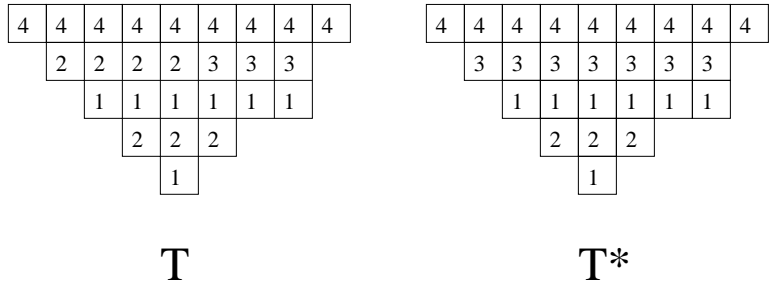}
\caption{Two minimal bar tableaux of shape 97631}
\label{fig:minbartableau}
\end{figure}

For example Figure \ref{fig:minbartableau} shows a minimal bar tableau
$T$ of shape 97631 and a minimal bar tableau $T^*$ of the same shape
with no bar boundaries an even number of squares along any row  (we will
verify later that these tableaux are minimal).

\begin{proof}
Let $T$ be a minimal bar tableau of $\lambda$.  In each row $r_k$ of
$T$, at the last bar boundary an even number of squares 
along a row, let $b$ be the bar which begins to the right of the
boundary.  Say that $b$ is labelled by $j$.  Relabel the squares to the left
of the boundary with $j$.  This preserves
the ordering on labels and the parity of $b$.  The partitions
remaining if we remove all squares labelled by integers larger than $i
(> j)$ will be the same as before and have distinct parts.  The
partitions remaining if we remove all squares labelled by $i (<j)$ will
not contain row $r_k$ but will otherwise have the same (distinct) parts
as before.
\end{proof}

\begin{lemma}\label{l:prepstruct}
Let $T^*$ be a minimal bar tableau of $\lambda$ such that there is no
bar boundary an even number of squares along any row.  Then if row $r_k$
is odd, it is labelled entirely by one label $j$.  If row $r_k$ is even,
it is labelled entirely by one label $j$ or it has exactly two labels
each occurring an odd number of times.
\end{lemma}

\begin{proof}
If row $r_k$ is odd and has more than one label, the second bar must be
of type 1. Therefore the second bar must be odd, forcing the first 
bar to be even which
is a contradiction.

If row $r_k$ is even and has more than two labels, the final two bars
are both of type 1 and so must be odd, forcing there to be a bar
boundary an even number of squares along the row, a contradiction.
\end{proof}

\section{Number of strips in a Minimal Tableau}

We use the results from the previous section to give a
count of how many bars are needed in a minimal bar tableau.
Define the \emph{shifted rank} of a shape $\lambda$, denoted
$\mathrm{srank}(\lambda)$,  to be the number of bars  
in a minimal bar tableau of $\lambda$.  Given an integer $a$, define
$a$ mod 2 to be 1 if $a$ is odd and 0 if $a$ is even.

\begin{theorem}\label{t:count}
Given a shape $\lambda$, let $o$ be the number of odd rows of
$\lambda$ and $e$ be the number of even rows.  Then $\mathrm{srank}
(\lambda) =
\mathrm{max}(o, e + ( \ell(\lambda) \ \mathrm{mod} \ 2))$.
\end{theorem}
For example, if $\lambda = 97631$, we have $o = 4, \ell(\lambda)=5$
and $e=1$.  So $\mathrm{srank}
(\lambda) = \mathrm{max}(4,1+1) = 4$ which verifies that the tableaux
shown in Figure \ref{fig:minbartableau} are indeed minimal.  If
$\lambda = 432$, we have $o = 1, \ell(\lambda)=3$
and $e=2$.  So $\mathrm{srank}
(\lambda) = \mathrm{max}(1,3) = 3$.  Such a tableau is illustrated below.
\begin{figure}[h]
\centering
\includegraphics[width=1.0in]{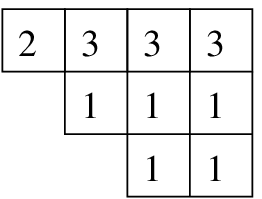}
\label{fig:minbartableau2}
\end{figure}
\begin{proof}
Let $T$ be a minimal bar tableau of $\lambda$.  Preprocess $T$ into
$T^*$ so that there are no bar boundaries an even number of squares
along any row.  This must preserve the number of bars.  Bars of type 3
consist of one even initial bar and one odd initial bar, and so by
Lemma \ref{l:prepstruct} must be an entire even row and an entire
odd row, or an entire even row and the initial odd bar of some other
even row.

First assume that $o \ge e$.  Note that if $o = e$, then $\ell(\lambda)
\mbox{ mod } 2 = 0$.  So when $o \ge e$,  $\mathrm{max}(o, e + ( \ell(\lambda)
 \ \mathrm{mod} \ 2))) = o$.  We claim that the bars of type 3 all
 consist of entire even row and entire odd row pairs, and that there
 are exactly $e$ of them.  

From the observations above, the number of bars of type 3
cannot be larger than $e$.  Suppose that there is a bar of type 3
consisting of an entire even row and the initial odd bar of some other
even row.  Since $o \ge e$, there must also be two other odd rows, not
parts of bars of type 3,
each labelled entirely by some label (by Lemma \ref{l:prepstruct}).
The total number of bars in these 4 rows is 4.  So if we relabel (with
new large labels) these four rows as two bars of type 3, we save two
bars, contradicting the minimality of $T^*$  (note that relabelling
entire rows with some integer larger than any current label preserves
legality, provided the parities are correct).  We illustrate this
(impossible) situation below, and show the more economical version.  Thus
there are no such 
bars of type 3. 
\begin{figure}[h]
\centering
\includegraphics[width=3.0in]{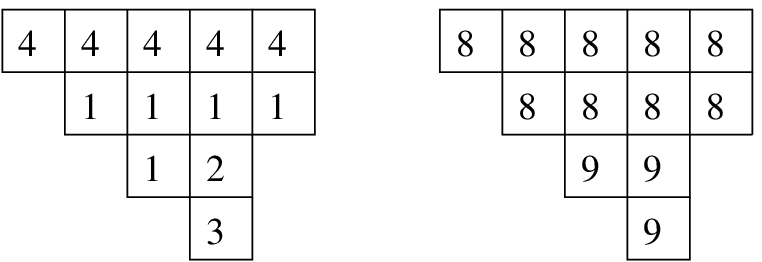}
\label{fig:minbartableau3}
\end{figure}

Now suppose that the number of bars of type 3 is
smaller than $e$; thus there is some even row $r_1$ of the tableau which is not
part of a bar of type 3.
Also there is an odd row $r_2$ which is not part of a type 3 bar
(since $ o \ge e$).    But we could relabel both these rows with some
new large label saving at least one bar  and contradict the minimality
of $T^*$. 

So there are exactly $e$ bars of type 3, filling $2e$
rows of $\lambda$. 
The remaining $o-e$ rows are odd and so must each be completely filled
by a unique label.  So the total number of bars is $o$ as
required.

\vskip 1.5em

Now assume that $e \ge o$. First we show that we can relabel 
so that every odd row is 
part of a bar of type 3.  So suppose $r_3$ is an odd row which is not part
of a bar of type 3.  

\emph{Claim.} There is an even row $r_4$, completely filled by a label,
which is part of a bar of 
type 3 with the initial part of some other even row $r_5$.  

\emph{Proof of claim.}  Assume by way of contradiction that there
is not, i.e. that the 
completely filled even rows are all parts of bars of type 3 with
complete odd
rows. But there must be at least one even row $r_4$ (since $e \ge o$
and $r_3$ is not part of a bar of type 3) which is not part of a bar
of type 3 with a complete odd row.  So $r_4$ must not be part of a bar
of type 3 at all (by our claim assumption).  
But then we could relabel $r_4$ and
$r_3$ entirely with some new large label
and save a bar, a contradiction.  This proves our claim.  
 
Relabel $r_4$ and $r_3$ with some new large label.  This
leaves an odd number of squares in the initial part of row $r_5$, and so
preserves legality.
These two rows are now a valid bar of type 3, and this process did not
cost us any bars.  We illustrate one step of this process below.
Simply iterate this process until there are no 
odd rows which are not part of bars of type 3.  This proves that we can relabel
so that every odd row is part of a bar of type 3.
\begin{figure}[h]
\centering
\includegraphics[width=2.5in]{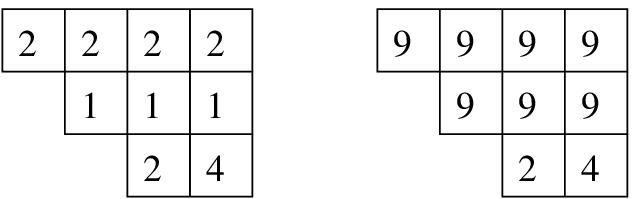}
\label{fig:minbartableau4}
\end{figure}

So every odd row is part of a bar of type 3, filling $2o$ rows of
$\lambda$.  All (except for possibly one) of the even rows remaining must be paired up
with another remaining even row, and each pair must contain one bar of type 3
(filling one entire row and the odd initial part of the other row) and
one bar of type 1 (filling the odd final part of the other row).
If they were not, we would have two even rows
costing 4 strips, and could reduce the number of strips by relabelling
as above with large new numbers.  The extra
row exists only when $\ell(\lambda)$ is odd, and costs two bars (i.e. one
extra).  This situation is illustrated on the right hand side of the
above figure.  So we
have $o + e - o + (\ell(\lambda) \mbox{ mod } 2)$ strips as required.
\end{proof}

It remains an open problem to count how many minimal shifted tableau
there are for a given shape $\lambda$.  Also, it would be natural to
generalise Theorem \ref{t:count} to the skew shifted case.  It is
possible to define skew bar tableaux, but the minimal number
of bars therein remains an open question (see \cite{Cl03} for more
details).

\section{Projective Representations of the Symmetric Group}

Here we recall some facts about the projective representations of the
symmetric group.  We follow the
treatment of Stembridge \cite{St88}.

A projective representation of a group $G$ on a vector space $V$ is a
map $P:G \rightarrow GL(V)$ such that
$$ P(x) P(y) = c_{x,y} P(xy) \ \ \ \  (x,y \in G)$$
for suitable (nonzero) scalars $c_{x,y}$.  For the symmetric group,
the associated Coxeter presentation shows that a representation $P$
amounts to a collection of linear transformations 
$\sigma_1, \ldots, \sigma_{n-1} \in GL(V)$ (representing the adjacent
  transpositions) such that $\sigma_j^2, (\sigma_j \sigma_{j+1})^3$,
  and $(\sigma_j \sigma_k)^2$ (for $ |j-k| \ge 2$) are all scalars.
  The possible scalars that arise in this fashion are limited.  Of
  course, one possibility is that the scalars are trivial; this occurs
  in any ordinary linear representation of $S_n$.  According to a
  result of Schur \cite{Sc11}, there is only one other possibility
  (occurring only when $n \ge 4$); namely
\begin{equation}
 \sigma_j^2 = -1; \ \ (\sigma_j \sigma_k)^2 = -1 \  ( \mbox{for}\   |j-k|
\ge 2); \ \  (\sigma_j \sigma_{j+1})^3 = -1.
\label{eq:sigmas}
\end{equation}
All other possibilities can be reduced to this case or the trivial
case by a change of scale.  See \cite{Jo89}, \cite{St89} for details.

It is convenient to regard $\sigma_1, \ldots, \sigma_{n-1}$ as
elements of an abstract group, and to take \ref{eq:sigmas} as a set of
defining relations.  More precisely, for $n \ge 1$ let us define
$\tilde{S}_n$ to be the group of order $2 \cdot n!$ generated by
$\sigma_1, \ldots, \sigma_{n-1}$ (and $-1$), subject to the relations
\ref{eq:sigmas}, along with the obvious relations $(-1)^2 = 1,
(-1)\sigma_j = \sigma_j(-1)$ which force $-1$ to be a central
involution.  By Schur's Lemma, an irreducible linear representation of
$\tilde{S}_n$ must represent $-1$ by either of the scalars $+1$ or
$-1$.  A representation of the former type is a linear representation
of $S_n$, whereas one of the latter type corresponds to a projective
representation of $S_n$ as in \ref{eq:sigmas}.   We will refer to any
representation of $\tilde{S}_n$ in which the group element $-1$ is
represented by the scalar $-1$ as a \emph{negative representation} of
$\tilde{S}_n$.

Next we review the characters of the irreducible negative representations of
$\tilde{S}_n$. Define $\calp(n)$ to be the set of all partitions of
$n$.  We say that a partition $\lambda$ is \emph{odd} if and
only if the 
number of even parts in $\lambda$ is odd, and is \emph{even} if and
only if it is not odd.  Thus, the parity of a permutation agrees with
the parity of its cycle type.  The parity of $\lambda$ is also the
parity of the integer $| \lambda | + \ell(\lambda)$.  Schur showed
that the irreducible negative representations are indexed by
partitions $\lambda$ with distinct parts.  Recall that if $P$ is an
irreducible negative representation indexed by $\lambda$ that the
character $\lang \lambda \rang$ is a class function $\lang \lambda
\rang:\tilde{S}_n \rightarrow \Q$ defined by $\lang \lambda
\rang (g) = \mathrm{trace}(P(g))$.

If $g=\pm \sigma_{i_1} \sigma_{i_2} \cdots$, let $\pi \in
\calp(n)$ be the cycle type (in $S_n$) of $\sigma_{i_1} \sigma_{i_2}
\cdots$.  In the sequel we will evaluate $\lang \lambda
\rang (\pi)$ instead of $\lang \lambda
\rang (g)$. Define $\calp ^0
(n) $ to be all partitions of $n$ such that all parts are odd.

\begin{theorem}[Schur 1911 \cite{Sc11}]
Let $\lambda \in \cald (n)$ have length $\ell$, and let $\pi \in
\calp(n)$.
\begin{enumerate}
\item  Suppose that $\lambda$ is odd.  If $\pi$ is neither in $\calp
  ^0 (n)$ nor equal to $\lambda$ 
  then $\lang \lambda \rang (\pi) = 0$.
\item  \label{it:odd} Suppose that $\lambda$ is odd.  If $\pi$ equals
  $\lambda$ then   
$$\lang \lambda \rang (\pi) = \pm
  i^{(n-\ell+1)/2 } (\lambda_1 \lambda_2 \cdots \lambda_{\ell} /2 )^ {1/2}.$$
\item Suppose that $\lambda$ is even.  If $\pi$ is not in $\calp ^0
  (n)$ then  $\lang \lambda \rang (\pi) = 0$. 
\end{enumerate}
\end{theorem}

For example we consider the situation when  $n=6$ and $\lambda=321$.   Then $\lang \lambda \rang (\pi) =
0$ when $\pi$ is $(6), (42), (411), (222), (2211)$ or
$(21111)$, as these partitions all have one even part.   
The second fact gives $\lang \lambda \rang (\pi) = \sqrt{3}$ when $\pi
= (321)$.  If
$\lambda=51$ then $\lang \lambda \rang (\pi) = 
0$ when $\pi$ is $(6),(42),  (411),  (321),  (222),
(2211)$ or 
$(21111)$. 

A combinatorial rule for calculating the characters not specified by
Schur's theorem was given by Morris; it is the projective analogue of
the Murnaghan-Nakayama rule.
\begin{theorem}[Morris 1962 \cite{Mo62}]
Let $\lambda \in \cald (n)$ have length $\ell$.
Suppose that $\pi \in \calp ^0
(n) $ and that $\pi$ contains $r$ at least once.  Define $\pi' \in \calp ^0
(n-r)$ by removing a copy of $r$ from $\pi$.  Then 
$$\lang \lambda \rang ( \pi ) = \sum_{i \in I(\lambda,r)} n_i \lang
\lambda (i,r) \rang (\pi ' ),$$
where
$$  n_i = \left\{
\begin{tabular}{ll}
$(-1)^{j-i} 2^{1-\varepsilon ( \lambda )} $ & $ \mbox{if } i \in I_+ ; $\\
$(-1)^{\ell -i} $  & if $ i \in I_0 $; \\
$(-1)^{j-i+\lambda_i} 2^{1-\varepsilon ( \lambda )} $ & $ \mbox{if } i
\in I_- . $\\ 
\end{tabular}
\right.
  $$
(The integer $j$ is that occurring in the definitions of $I_{\pm}$,
and $\varepsilon(\lambda)$ is the parity of $\lambda$; i.e. 0 or 1.)
\end{theorem}

For example if $n=6$, $\lambda=(51)$ and $r=1$, we have
$\varepsilon(\lambda)=0, I_+ = \{1\}, I_0 = \{2\} $ and $  I_- =\emptyset
.$  So $I(\lambda,r) = \{1,2\}$, and we have
\begin{eqnarray*}
 \lang 51 \rang ( 1^6 ) &=& (-1)^{1-1} 2^{1-0}\lang
41 \rang (1^5 ) + (-1)^{2-2}  \lang
5 \rang (1^5 )\\
&=& 2\lang 41 \rang (1^5 ) +  \lang 5 \rang (1^5 ).
\end{eqnarray*}

Expand this sum into a sum over all possible bar tableaux.
Define the \emph{weight} of a tableau wt$(T)$ to be the product 
of all the powers of $-1$ and $2$  which appear.  Then we have 
$$\lang \lambda \rang ( \pi ) = \sum_{T} wt(T),$$
summed over all bar tableaux of shape $\lambda$ and type $\pi$.
We know that the shifted rank of $\lambda$ is the minimum number of bars
needed in a bar tableau of shape $\lambda$. So we obtain the
following result as a corollary to Theorem \ref{t:count}:

\begin{corollary}
Given a shape $\lambda$ of shifted rank $k$ and a shape $\pi$ such that
$\ell(\pi) < k$, we have $\lang \lambda \rang (\pi) = 0$.   
$\Box$\end{corollary}

\section{Schur Q-Functions}\label{subsec:qfunctions}

We begin with Schur's original inductive definition of the
$Q_{\lambda}$ functions. 
Denote the monomial symmetric functions by $m_{\lambda}$
and define symmetric
functions 
$q_k$ of
degree $k$ by 
$$ q_k = \sum_{\lambda \in \calp(k)} 2^{\ell(\lambda)} m_{\lambda}.$$
Now we can state the base cases for the inductive definition.  Put
$Q_{(a)} = q_a$ and
$$Q_{(a,b)} = q_a q_b + 2 \sum_{n > 0} (-1)^n q_{a+n}q_{b-n}.$$
Inductively we define 
$$ Q_{\lambda_1, \ldots , \lambda_{2k+1}} = \sum_{i=1}^{2k+1}
(-1)^{i+1} q_{\lambda_i} Q_{\lambda_1, \ldots, \hat{\lambda}_i, \ldots,
  \lambda_{2k+1}} $$
and
$$ Q_{\lambda_1, \ldots , \lambda_{2k}} = \sum_{i=2}^{2k}
(-1)^{i} Q_{\lambda_1,\lambda_i} Q_{\lambda_2, \ldots, \hat{\lambda}_i, \ldots,
  \lambda_{2k}}.$$

The $Q_{\lambda}$ may also be defined as the specialization at $t=-1$
of one of the two equivalent defining formulae for Hall-Littlewood
polynomials; see \cite[III (2.1) (2.2)]{Ma95}.  Let $S_r$ act on
$X=\{x_1, \ldots, x_r \}$ by permuting the variables, so that, when
$\ell \le r$, the Young subgroup $S_1^{\ell} \times S_{r-l}$ fixes
each of $x_1, \ldots, x_{\ell}$.  Let $\lambda$ be a strict partition of
length $\ell$.  If $\ell \le r$, then
$$Q_{\lambda}(x_1 , \ldots , x_r) = 2^{\ell} \sum_{ [w] \in S_r /
  S_1^{\ell} \times S_{r-\ell}} w \left \{ x_1^{\lambda_1} \cdots
x_{\ell}^{\lambda_{\ell}}  \prod_{i=1}^{\ell} \prod_{j=i+1}^{r}
\frac{x_i+x_j}{x_i-x_j} \right \}.$$
If $\lambda$ has length greater than $r$, then $Q_{\lambda}(x_1,
\ldots, x_r)=0$.    The $Q_{\lambda}$ symmetric functions are obtained
by taking the limit as the number of variables becomes infinite (for a
mathematically precise definition of this limit see \cite{Ma95}).

Schur \cite{Sc11} defined these Q-functions in order to
study the projective representations of symmetric groups.  The
fundamental connection is given by the following theorem.  Let
$m_i(\lambda)  = \# \{j: \lambda_j = i \}$, the
number of parts of $\lambda$ equal to $i$.  Define $z_{\lambda} =
1^{m_1(\lambda)} m_1(\lambda)! 2^{m_2(\lambda)} 
  m_2(\lambda)!\cdots$.
  Denote the power sum
    symmetric functions 
  by $p_{\lambda}$. 

\begin{theorem}[Schur, 1911]
$$Q_{\lambda}  = \sum_{\pi \in \calp ^ 0 (n)}
2^{[\ell(\lambda)+\ell(\pi) + \varepsilon(\lambda)]/2} \lang
\lambda \rang (\pi) \frac{p_{\pi}}{z_{\pi}}.$$
\end{theorem}

Again consider the example with $n=6$ and $\lambda=(51)$.  We have
\begin{eqnarray*}
Q_{51} &=& 2^{[2+6+0]/2} \lang 51 \rang (1^6)  \frac{p_{1^6}}{z_{1^6}}+
 2^{[2+4+0]/2} \lang 51 \rang (1^3 3)  \frac{p_{1^3 3}}{z_{1^3 3}}+ \\
& & 2^{[2+2+0]/2} \lang 51 \rang (15)  \frac{p_{1 5}}{z_{15}}+
2^{[2+2+0]/2} \lang 51 \rang (3^2)  \frac{p_{3^2}}{z_{3^2}} \\
&=& 2^{4} 16 \frac{p_{1^6}}{6!}+
2^3 2 \frac{p_{1^3 3}}{3! 3} -
2^2 1 \frac{p_{1 5}}{5}-
2^2 2 \frac{p_{3^2}}{18}\\
&=& \frac{16}{45} p_{1^6} +  \frac{8}{9}p_{1^3 3}
-\frac{4}{5} p_{1 5} -  \frac{4}{9}p_{3^2}
\end{eqnarray*}

Define deg($p_i$) $=1$, so deg($p_{\nu}) = \ell(\nu)$.  Then Theorem
\ref{t:count} gives us the following corollary.
\begin{corollary}
The terms of lowest degree in $Q_{\lambda}$ have degree at least
$\mathrm{ srank 
  }  (\lambda)$. $\Box$
\end{corollary}

In our example srank$(51) = 2$ and the $p_{\nu}$ satisfy $\ell(\nu)
\ge 2$.  
Equivalently we can examine a specialization of the principal
specialization of $Q_{\lambda}$, i.e.
$$ps^1_t(Q_{\lambda})  = Q_{\lambda}(\underbrace{1,1,\ldots,1}_{t \  1's})=
Q_{\lambda}(1^t).$$ 
Since $p_{\nu}(1^t) =
t^{\ell(\nu)}$, we can rephrase the above result.
\begin{corollary}
$Q_{\lambda}(1^t)$ is divisible by
$t^{ \, \mathrm{srank} (\lambda)}$. $\Box$
\end{corollary}

The following conjecture has been computationally verified (using
John Stembridge's SF Package for Maple \cite{St01}) for all
partitions $\lambda \vdash n$ for $1\le n \le 12$.
\begin{conjecture}
The terms of lowest degree in $Q_{\lambda}$ have degree exactly
$\mathrm{ srank 
  }  (\lambda)$.
\end{conjecture}

\end{document}